\documentclass{amsart}          

\usepackage{amsmath,amsthm}     
\usepackage{amssymb}            
\usepackage{euscript}           
\usepackage{enumerate,calc}
\usepackage[matrix,arrow,curve,frame]{xy}    

\xymatrixcolsep{1.9pc}                          
\xymatrixrowsep{1.9pc}
\newdir{ >}{{}*!/-5pt/\dir{>}}                  

\def\longfib{\DOTSB\relbar\joinrel\twoheadrightarrow}

\newtheorem{thm}[subsection]{Theorem}
\newtheorem{defn}[subsection]{Definition}
\newtheorem{prop}[subsection]{Proposition}

\newtheorem{cor}[subsection]{Corollary}
\newtheorem{lemma}[subsection]{Lemma}

\theoremstyle{definition}  

\newtheorem{remark}[subsection]{Remark}

  {\end{list}}

\newcommand{\dfn}{\textbf} 

\newcommand{\mdfn}[1]{\dfn{\mathversion{bold}#1}} 

\newcommand{\tens}              {\otimes}               
\newcommand{\iso}               {\cong}  

\newcommand{\cat}{\EuScript}    
\newcommand{\cA}{{\cat A}}      
\newcommand{\cC}{{\cat C}}
\newcommand{\cD}{{\cat D}}
\newcommand{\cE}{{\cat E}}
\newcommand{\cF}{{\cat F}}

\newcommand{\cM}{{\cat M}}
\newcommand{\cN}{{\cat N}}

\newcommand{\cR}{{\cat R}}

\newcommand{\Top}{{\cat Top}}

\newcommand{\sSet}{s{\cat Set}}

\newcommand{\Pre}{Pre}

\newcommand{\ho}{\text{Ho}\,}

\newcommand{\field}[1]  {\mathbb #1} 

\newcommand{\R}         {\field R}
\newcommand{\LL}        {\field L}

\newcommand{\Z}         {\field Z}

\DeclareMathOperator*{\colim}{colim}
\DeclareMathOperator*{\hocolim}{hocolim}
\DeclareMathOperator*{\badhocolim}{badhocolim}

\DeclareMathOperator{\coeq}{coeq}

\newcommand{\ra}{\rightarrow}                   
\newcommand{\lra}{\longrightarrow}              
\newcommand{\lla}{\longleftarrow}               
\newcommand{\llra}[1]{\stackrel{#1}{\lra}}      
\newcommand{\llla}[1]{\stackrel{#1}{\lla}}      

\newcommand{\we}{\llra{\sim}}                   
\newcommand{\bwe}{\llla{\sim}}
\newcommand{\cof}{\rightarrowtail}              
\newcommand{\fib}{\twoheadrightarrow}           
\newcommand{\trfib}{\stackrel{\sim}{\longfib}}
\newcommand{\trcof}{\stackrel{\sim}{\cof}}

\newcommand{\inc}{\hookrightarrow}              
\newcommand{\dbra}{\rightrightarrows}           

\newcommand{\Id}{Id}                            

\newcommand{\norm}[1]{\mid \! #1 \! \mid}       
\newcommand{\mM}{\underline{\cM}}

\newcommand{\scop}{{\Delta^{op}}}
\newcommand{\re}{Re}
\newcommand{\Rea}{\re\,}

\newcommand{\sing}{Sing}
\newcommand{\Sing}{\sing\,}

\newcommand{\assign}{\mapsto}

\newcommand{\ovcat}{\downarrow}

\newcommand{\adjoint}{\rightleftarrows}

\newcommand{\del}[1]{\Delta^{#1}}

\newcommand{\pt}{pt}

\newcommand{\ev}{ev_0}

\newcommand{\ovch}{(\cC\times\Delta \downarrow X)}
\newcommand{\ovcf}{(\cC\times\Delta \downarrow F)}
\newcommand{\CR}{\cC R}
\newcommand{\Cl}{\cC}

\newcommand{\fix}{\mbox{}\par\noindent}

\numberwithin{equation}{subsection}

\begin{document}

\title{Combinatorial model categories have presentations}

\author{Daniel Dugger}
\address{Department of Mathematics\\ Purdue University\\ West
Lafayette, IN 47907 } 

\email{ddugger@math.purdue.edu}

\begin{abstract}
We show that every combinatorial model category is Quillen equivalent
to a localization of a diagram category (where `diagram category'
means diagrams of simplicial sets).  This says that every
combinatorial model category can be built from a category of
`generators' and a set of `relations'.
\end{abstract}

\maketitle


\section{Introduction}
In the companion paper \cite{D2} we introduced a technique for
constructing model categories via generators and relations.  The two
main points were as follows:

\begin{enumerate}[(1)]
\item
From a small category $\cC$ one can construct a model category $U\cC$
which is in some sense the free model category generated by $\cC$.
This $U\cC$ is simply the category of diagrams $\sSet^{\cC^{op}}$ with
an appropriate model structure.  
\item Given a set of maps $S$ in $U\cC$, one can form the
\dfn{localization} $U\cC/S$---this is the closest model category to
$U\cC$ in which the maps from $S$ have been added to the weak
equivalences.  We regard this process of localization as `imposing
relations' into the model category $U\cC$.  
\end{enumerate}

When a model category can be built from generators and relations we
say that it has a \dfn{small presentation}.  More precisely, a
presentation for a model category $\cM$ consists of (1) a small
category $\cC$, (2) a set of maps $S$ in the diagram category $U\cC$,
and (3) a specified Quillen equivalence $L:U\cC/S \adjoint \cM:R$.
As in \cite{D2} we will denote a Quillen pair as a map of model
categories in the direction of the left adjoint, so that our
presentation takes on the form $U\cC/S \we \cM$.  

There are certainly model categories which cannot be given
presentations, but the majority of those one encounters 
can indeed be built up in this way.  In this paper we will deal with a
very broad class called the {\it combinatorial\/} model categories,
which were introduced by Jeff Smith.  These include essentially any
model category of algebraic origin, as well as anything constructed in
some way out of simplicial sets.  Our aim is to prove the following
result, announced in \cite{D2}:

\begin{thm}
\label{th:main}
Every combinatorial model category has a small presentation.  
\end{thm}

This has the following corollary:

\begin{cor}
\label{co:main}
Every combinatorial model category is Quillen equivalent to one which
is simplicial, left proper, and (this is slightly harder) in which
every object is cofibrant.
\end{cor}

In \cite{D1} it was proven, using very different methods, that every
left proper, combinatorial model category is Quillen equivalent to a
simplicial model category.  The above corollary offers a slight
improvement on this, in that it eliminates the left-properness
assumption.

We close this introduction with a word about the proof of
Theorem~\ref{th:main}.  To establish some intuition, consider the
problem of giving a presentation for an abelian group $A$.  The first
thing one does is to find a surjection $\Z^r \fib A$, and if $R$
denotes the kernel then one automatically has the presentation $\Z^r/R
\iso A$.  Following this analogy our approach will be to define what
it means for a map of model categories $\cM\ra \cN$ to be a
`surjection' (\ref{de:surj}), and we'll see that once one has a
surjection $U\cC \ra \cM$ then a presentation follows automatically
(\ref{pr:quot}).  So getting the surjection will be the tricky part,
and this depends on carefully choosing the category $\cC$.

Since we think of $\cC$ as a category of `generators', it's natural to
try and choose it as a subcategory of $\cM$.  Intuitively, $\cC$
should be big enough so that every object $X\in \cM$ can be built up
as a certain homotopy colimit of objects from $\cC$, but this turns
out to be more delicate than it sounds (see section 4 for a precise
statement).  In the case when $\cM$ is a {\it simplicial\/}
model category, though, this plan can indeed be carried out.

When $\cM$ is not simplicial we have to be more clever.  Instead of
choosing $\cC$ as a subcategory of $\cM$, we are forced to fatten
$\cM$ a little by looking at the cosimplicial objects $c\cM$.  By
choosing an appropriate subcategory $\cC$ of $c\cM$ (section 6) we are
able to get our surjective map $U\cC \ra \cM$.

\subsection{Overview}

Section 2 contains a brief discussion of combinatorial model
categories, especially several key properties that will appear
throughout the paper.  In section 3 we give a definition of
`homotopically surjective' maps, and we prove both
Theorem~\ref{th:main} and its Corollary assuming the existence of a
surjective map $U\cC \ra \cM$.  The remainder of the paper is the
quest for this surjective map.

Section 4 deals with `canonical homotopy colimits'.  Given a map
$\gamma\colon \cC \ra \cM$ and a fibrant object $X\in \cM$ the
canonical homotopy colimit \mdfn{$\hocolim \ovch$} is an object built
out of all the information contained in the homotopy function
complexes $\mM(\gamma c,X)$, for all $c\in \cC$---it is a kind of
homotopical approximation to $X$ based on $\cC$.  We show
(\ref{co:hocosam}) that finding a surjective map $U\cC \ra \cM$ is
equivalent to finding a functor $\gamma\colon \cC \ra \cM$ such that
these homotopical approximations to $X$ always give back $X$ itself
(up to weak equivalence, of course).

In section 5 we produce the required surjective map in the case when
$\cM$ is a {\it simplicial\/} model category.  This case is fairly
simple based on our work so far.  Section 6 handles the more general
case---the ideas are similar to those from section 5, but with an
extra level of complication.

Sections 7 and 8 contain some of the auxiliary proofs postponed from
previous sections.  Finally, since much of the paper is spent working
with various homotopy colimits, we have for convenience enclosed an
appendix recalling some of the basic properties we need.  In
particular, there are several instances in the paper where we have to
identify two homotopy colimits over different indexing categories, and
the key result letting us do this is Proposition~\ref{pr:hocored}.

\subsection{Notation}
This paper is intended as a companion to \cite{D2}, and we assume a
general familiarity with the notation and results of sections 2, 3, and
5 of that paper.  We deal quite a bit with overcategories here, so
recall that if $F\colon \cC \ra \cD$ is a functor and $X\in \cD$ then
$(F\ovcat X)$---often written $(\cC\ovcat X)$ by abuse---is the
category whose objects are pairs $[c, Fc\ra X]$ where $c\in \cC$ and
$Fc \ra X$ is a map in $\cD$.  The morphisms of $(\cC\ovcat X)$ are
the obvious candidates.

\subsection{Acknowledgements}
I would like to express my thanks to Jeff Smith for sharing his
work on combinatorial model categories with me, as well as for several
useful conversations about the results in this paper.

\section{Combinatorial model categories}

In this section we review the definition of combinatorial model
categories, due to Jeff Smith, together with several of their important
properties.  The main theorem of this paper (\ref{th:main}) is the
homotopy-theoretic analog of a standard result about locally
presentable categories, recalled in (\ref{se:locdia}).

\medskip

\begin{defn}
A model category $\cM$ is called \dfn{combinatorial} if it is
cofibrantly-generated and the underlying category is locally
presentable.
\end{defn}

This definition is surprisingly powerful considering how simple it is.
We'd better recall what all the words mean, though.  The notion of
{\it cofibrantly-generated\/} model category is standard by now, and may
be found in \cite[Def. 2.1.17]{Ho}---it requires that there are basic
sets of cofibrations and trivial cofibrations which one can use to do
the small object argument.  The notion of a category being {\it
locally presentable\/} is less familiar to homotopy-theorists, so here's
the definition:

\begin{defn}
A category $\cC$ is \dfn{locally presentable} if it is co-complete,
and if there is a regular cardinal $\lambda$ and a set of objects
$\cA$ in $\cM$ such that
\begin{enumerate}[(i)]
\item Every object in $\cA$ is small with respect to $\lambda$-filtered
colimits, and
\item Every object of $\cM$ can be expressed as a $\lambda$-filtered
colimit of elements of $\cA$. 
\end{enumerate} 
\end{defn}

For background on locally presentable categories one may consult
\cite[Section 1.B]{AR} or \cite{B}.  The condition that an object be
small with respect to $\lambda$-filtered colimits is called
\mdfn{$\lambda$-presentable} in \cite{AR}, but we will follow Smith
and call it \mdfn{$\lambda$-small}.  Locally presentable categories
have the following important characteristics:

\begin{enumerate}[(1)]
\item For every object $A$, there exists a regular cardinal $\lambda$
such that $A$ is $\lambda$-small.
\item For each regular cardinal $\lambda$, the $\lambda$-small objects
in $\cC$ have a {\it set\/} of representatives with respect to
isomorphism---we'll use \mdfn{$\cC_\lambda$} to denote the full
subcategory determined by any such set.
\end{enumerate}

The following proposition brings together the properties of
combinatorial model categories we will need in this paper.  Most of
these statements are due to Smith, and should one day appear in
\cite{Sm}.  For the reader's convenience we provide proofs (or
sketches of proofs, when we are lazy) in section 7.

\begin{prop}\label{li:combmc}\label{pr:combpr}
Let $\cM$ be a combinatorial model category.
\begin{enumerate}[(i)]
\item There exist cofibrant- and fibrant-replacement functors which
preserves sufficiently large filtered colimits.
\item Sufficiently large filtered colimits of weak equivalences are
again weak equivalences: if $\lambda$ is a sufficiently large regular
cardinal, $I$ is a $\lambda$-filtered indexing category, and if
$D_1,D_2\colon I \ra \cC$ are diagrams with a natural
weak equivalence $D_1 \ra D_2$, then $\colim D_1 \ra \colim D_2$ is
also a weak equivalence.
\item There exist functorial factorizations of maps $X\ra Y$ as $X
\trcof \tilde{X} \fib Y$ and $X \cof \tilde{Y} \trfib Y$ with the
following property: for sufficiently large regular cardinals $\mu$, if
$X\ra Y$ is a map between $\mu$-small objects then both $\tilde{X}$
and $\tilde{Y}$ are $\mu$-small as well.  
\end{enumerate}
\end{prop}

Notice that the above properties are well-known for the model category
of simplicial sets.  They in some sense say that for a combinatorial
model category the interesting part of the homotopy theory is all
concentrated within some small subcategory--- beyond sufficiently
large cardinals the homotopy theory is somehow `formal'.  Model
categories of the form $U\cC/S$ certainly have this property (as they
are combinatorial), and this observation explains why not every model
category can have a small presentation.

\subsection{Locally presentable categories and diagram categories}
\label{se:locdia}
The proof of Theorem~\ref{th:main} is somewhat involved, and it will
help for us to establish a little background.  Like many of the
results in \cite{D2}, the theorem is once again a homotopy-theoretic
analog of a standard result in category theory.

Suppose that $\cC$ is a category and $F\colon\cA\ra \cC$ is a functor
(in many cases this will be the inclusion of a subcategory).  For any
$x\in \cC$ consider the overcategory $(\cA\downarrow x)$, together
with the canonical diagram $(\cA\downarrow x)\ra \cC$ which sends $[a,
Fa\ra x]$ to $Fa$.  The colimit of this diagram (when it exists) is
called the \mdfn{canonical colimit of $x$ with respect to $\cA$} and
we'll denote it \mdfn{$\colim (\cA\ovcat x)$}.

A locally presentable category $\cC$ has the following important
property: for large enough regular cardinals $\lambda$, every $x$ in
$\cC$ is isomorphic to its canonical colimit with respect to the
subcategory $\cC_\lambda$---one says that $\cC_\lambda$ is
\dfn{dense} in $\cC$.

When $\cC$ is locally presentable it is co-complete, and so the
inclusion $\cC_\lambda \inc \cC$ extends to an adjoint pair
\[ \re:\Pre(\cC_\lambda) \adjoint \cC:\sing \]
as in \cite[Prop. 2.1]{D2}.  It's not hard to check that
$\re(\sing\, x)$ is precisely the canonical colimit of $x$ with
respect to $\cC_\lambda$, and so the map $\re(\sing\, x) \ra x$ is an
isomorphism (again, for large enough $\lambda$).  A standard result in
the theory of locally presentable categories roughly says that there
is a `localization functor' $L\colon\Pre(\cC_\lambda) \ra
\Pre(\cC_\lambda)$ such that the above adjoint pair restricts to the
image of $L$ and becomes an equivalence of categories---in other
words, $\cC$ is equivalent to a full, reflective subcategory of the
diagram category $\Pre(\cC_\lambda)$ \cite[Prop. 1.46]{AR}.

Theorem~\ref{th:main} above is a direct homotopy-theoretic analog of
this result.  In section 4 we define the notion of canonical {\it
homotopy\/} colimits, and we will see in section 5 that combinatorial
model categories $\cM$ which are simplicial have the following
property: for sufficiently large regular cardinals $\lambda$ the
subcategory $\cM_\lambda$ is `homotopically dense', in that every $x$
in $\cM$ is weakly equivalent to its canonical homotopy colimit with
respect to $\cM_\lambda$.  This condition will allow us to get a
Quillen equivalence $U(\cM_\lambda)/S \adjoint \cM$ (actually we will
replace $\cM_\lambda$ by its subcategory of cofibrant objects, for
technical convenience, but this is not crucial).

For combinatorial model categories which are not simplicial the story
is slightly more complex, but the above ideas are still the central
points.  We refer the reader to the discussion which begins section 6
for more about this.

\section{Homotopically surjective maps}
In this section we will define what it means for a map of model
categories $\cM\ra \cN$ to be `surjective' (\ref{de:surj}), and we'll
see that a surjection from a universal model category automatically
yields a presentation (\ref{pr:quot}).  At the end of the section we
show how Theorem~\ref{th:main} and Corollary~\ref{co:main} will follow
as soon as one has found a surjective map of the form $U\cC \ra \cM$.

\begin{defn}
\label{de:surj}
A map of model categories $\cM\llra{L} \cN$ is \dfn{homotopically
surjective} if it has the following property: for every fibrant object
$X$ in $\cN$, and every cofibrant replacement $[RX]^{cof} \we RX$ for
$RX$, the induced map $L([RX]^{cof})\ra X$ is a weak equivalence in
$\cN$.  We often omit the word `homotopically' for brevity.
\end{defn}

Equivalently, the definition says that on the level of homotopy
categories the derived functors $\LL:\ho \cM \adjoint \ho \cN: \R$
are such that $\LL\circ \R$ is naturally isomorphic to the identity.

The following result says that for combinatorial model categories any
homotopically surjective map may be localized so as to become a
Quillen equivalence.  (Note: The left-properness assumption on $\cM$
is there so that we may form the localization $\cM/S$, otherwise it is
unimportant.)  

\begin{prop} 
\label{pr:quot}
Let $\cM$ and $\cN$ be combinatorial model categories, where $\cM$ is
left proper.  Suppose that $L\colon\cM\ra \cN$ is a surjective map.
Then there is a set of maps $S$ in $\cM$ which become weak
equivalences under $L^{cof}$, and such that the induced map $\cM/S \ra
\cN$ is a Quillen equivalence.
\end{prop}

\noindent
(Recall that $L^{cof}$---which we call the left-derived functor of
$L$---denotes the result of pre-composing $L$ with some
cofibrant-replacement functor in $\cM$.)

\begin{proof}
Choose a regular cardinal $\lambda$ which is large enough so that the
following are true:
\begin{enumerate}[(1)]
\item $\cM_\lambda$ is dense in $\cM$,
\item $\lambda$-filtered colimits of weak equivalences in $\cM$ are
again weak equivalences,
\item $\cM$ has a cofibrant replacement functor $A\ra A^{cof}$ which
preserves $\lambda$-filtered colimits,
\item $\cN$ has a fibrant replacement functor $X\ra X^{fib}$ which preserves
$\lambda$-filtered colimits,
\item the right adjoint $R$ to $L$ preserves $\lambda$-filtered
colimits (see \cite[Prop. 1.66]{AR}).
\end{enumerate}
Let $S$ be the set consisting of all the natural maps 
\[ A^{cof} \ra R([LA^{cof}]^{fib}) \] 
where $A\in \cM_\lambda$.  The condition that $\cM\ra \cN$ is
homotopically surjective shows that the derived functor of $L$ takes
maps in $S$ to weak equivalences, and so $L$ descends to a map $\cM/S
\ra \cN$.  It is readily checked that this new map is also
homotopically surjective.  To check that this is a Quillen equivalence
one must verify that for every object $X$ in $\cM$, the composite
$X^{cof}\ra R([LX^{cof}]^{fib})$ is a weak equivalence in $\cM/S$.
But any $X$ in $\cM$ is a $\lambda$-filtered colimit of objects in
$\cM_\lambda$ by assumption (1), and all the functors in sight commute
with such colimits by assumptions (3)--(5).  So the map in question is
a $\lambda$-filtered colimit of maps in $S$, which are weak
equivalences in $\cM/S$.  Finally, assumption (2) says that
$\lambda$-filtered colimits preserve weak equivalences in $\cM$, and
it's easy to check that this property is inherited by any localization
of $\cM$.  This completes the proof.
\end{proof}

The following proposition will be our focus in the rest of the paper.
Granting it for the moment, we can prove Theorem~\ref{th:main} and its
Corollary.

\begin{prop}
\label{pr:hosurj}
If $\cM$ is a combinatorial model category then there exists a small
category $\cC$ and a homotopically surjective map $U\cC \ra \cM$.
\end{prop}

\begin{proof}[Proof of Theorem~\ref{th:main}]
Let $\cC$ be the category guaranteed by the above proposition.  
The model category $U\cC$ is left proper and combinatorial, so by
Proposition~\ref{pr:quot} we can find a set of maps $S$ in $U\cC$
which become weak equivalences in $\cM$, and such that $U\cC/S \ra
\cM$ is a Quillen equivalence.
\end{proof}

\begin{proof}[Proof of Corollary~\ref{co:main}]
If $\cM$ is a combinatorial model category then the Theorem gives us a
Quillen equivalence $U\cC/S \we \cM$ for some $\cC$ and $S$.  
The point is that the universal model category $U\cC$ is simplicial
and left proper, and these properties are inherited by the
localization $U\cC/S$.

We must work a little harder to show that $\cM$ is Quillen equivalent
to a model category in which every object is cofibrant.  Recall that
the diagram category $\sSet^{\cC^{op}}$ has a \dfn{Heller model
structure} \cite{He} in which a map $D\ra E$ is a weak equivalence
(resp. cofibration) if $D(c) \ra E(c)$ is a weak equivalence
(resp. cofibration) for every $c\in \cC$.  The Heller model structure
is related to $U\cC$ by a Quillen equivalence $U\cC \we
\sSet^{\cC^{op}}_{H}$ (where the `H' is for `Heller').
This map will still be a Quillen equivalence when we localize,
so that we get a zig-zag of Quillen equivalences
\[ \sSet^{\cC^{op}}_{H}/S \bwe U\cC/S \we \cM.\]
But now the point is that the Heller model structure is simplicial,
left proper, and has the property that every object is
cofibrant; these properties all pass to the localization.
\end{proof}

The application of replacing combinatorial model categories by ones in
which everything is cofibrant was suggested to me by Jeff Smith.


\section{Canonical homotopy colimits}

In this section we introduce a homotopical generalization of canonical
colimits, which were discussed in (\ref{se:locdia}).  When $\cC$ is a
subcategory of $\cM$ then the {\it canonical homotopy colimit\/} of a
fibrant object $X$ with respect to $\cC$ is a certain `approximation'
to $X$ based on $\cC$: one takes all the information from the homotopy
function complexes $\mM(c,X)$ as $c\in\cC$ varies, and from this data
constructs the canonical homotopy colimit.  (In the general case $X$
need not be fibrant, and $\cC$ need not be a subcategory.)  The
importance for us is Corollary~\ref{co:hocosam}, which says that a map
$U\cC \ra \cM$ is homotopically surjective precisely if taking
canonical homotopy colimits with respect to $\cC$ always gives back
the original object up to weak equivalence.

\bigskip

Consider a functor $\gamma\colon\cC \ra \cM$ together with a
cosimplicial resolution $\Gamma\colon\cC \ra c\cM$ (see
\cite[Def. 3.2]{D2}).  If $c\in\cC$ then we'll use $\Gamma^n c$ to
denote the component of $\Gamma(c)$ lying in dimension $n$.  The
cosimplicial resolution induces a functor $\cC \times \Delta \ra \cM$
sending $(c,[n])$ to $\Gamma^n c$.  For each $X$ in $\cM$ one can form
the over-category $(\cC\times \Delta \downarrow X)$, together with the
canonical functor $(\cC\times\Delta \downarrow X)\ra \cM$.

\begin{defn}
The homotopy colimit of this functor is called the \dfn{canonical
homotopy colimit} of $X$ with respect to $\Gamma$ (or with respect to
$\cC$, if we are lazy), and it will be denoted \mdfn{$\hocolim
\ovch$}.
\end{defn}

As usual, talking about `the' canonical homotopy colimit is somewhat
inappropriate since the actual object depends on the framing used in
calculating the homotopy colimit---we will ignore this point of
etiquette, though.  Note that there are canonical maps
\[ \hocolim \ovch \ra \colim \ovch \ra X.\]
We will be very interested in the composite.

Everyone knows at least one example of a canonical homotopy colimit:
For the model category $\Top$, consider the
inclusion of the one-point category $\pt \inc \Top$ and its standard
cosimplicial resolution given by the topological simplices $\del{n}$.
The canonical homotopy colimit of a space $X$ with respect to this
subcategory turns out to be the same as the realization of the
singular complex of $X$.  

In general, $\hocolim \ovch$ is a kind of `homotopical approximation'
to $X$ based on the functor $\gamma\colon\cC \ra \cM$.  We look at all
ways of mapping $n$-fold homotopies $\Gamma^n c$ into $X$, and from
this data we concoct some strange object which is like a phantom image
of $X$ as seen through the eyes of $\gamma$.  In the above example
from topological spaces the natural maps $\hocolim \ovch \ra X$ are
all weak equivalences, but this will typically be far from true.

\medskip

From the above definition it is not immediately clear to what extent
the homotopy type of $\hocolim \ovch$ depends on the cosimplicial
resolution $\Gamma$, which contributed to the indexing category
$\ovch$.  We will show in a moment (\ref{co:canhoco}i) that if $X$ is
fibrant then choosing a different cosimplicial resolution yields a
weakly equivalent object.  We will also show (\ref{co:canhoco}ii) that
if $X \ra Y$ is a weak equivalence between fibrant objects, then the
induced map of canonical homotopy colimits is again a weak equivalence.
The key to proving these statements is the following result, which
says that canonical homotopy colimits can always be interpreted as
certain realizations of singular complexes (compare (\ref{se:locdia})):

\begin{prop}
Let $\re:U\cC \adjoint \cM:\sing$ be the Quillen pair induced by
$\Gamma$.  Then $\Rea^{cof} \Sing X$ is weakly equivalent to $\hocolim
\ovch$.
\end{prop}

\begin{proof}
Consider the Yoneda embedding $\cC \inc U\cC$, together with its
canonical cosimplicial resolution induced by the simplicial structure
on $U\cC$.  In \cite[Prop. 2.9]{D2} we showed that for any $F\in U\cC$
the natural map $\hocolim \ovcf \ra F$ gives a cofibrant-approximation
to $F$.  We aim to apply this in the case where $F$ is $\Sing X$.

It's easy to see using adjointness that the overcategory $(\cC\times
\Delta \ovcat \Sing X)$ is isomorphic to the overcategory $(\cC\times
\Delta \ovcat X)$.  It's then clear that applying the realization
$\re$ to $\hocolim (\cC\times\Delta \ovcat \Sing X)$ gives precisely
$\hocolim \ovch$.

So we've recovered $\hocolim \ovch$ by starting with $\Sing X$, taking a
certain cofibrant-approximation in $U\cC$, and then applying the
realization $\re$.  This is precisely what we needed to prove.
\end{proof}

\begin{cor}\label{co:canhoco}\fix
\begin{enumerate}[(i)]
\item If $X \ra Y$ is a weak equivalence between fibrant objects then
the induced map $\hocolim \ovch \ra \hocolim (\cC\times\Delta \ovcat
Y)$ is also a weak equivalence.
\item Suppose that $\Gamma'$ is another resolution for $\gamma$.  Then
the canonical homotopy colimits $\hocolim_{\Gamma} \ovch$ and
$\hocolim_{\Gamma'}  \ovch$ are weakly equivalent.
\end{enumerate}
\end{cor}

\begin{proof}
Part (i) follows directly from the fact that $(\re,\sing)$ is a
Quillen pair: the weak equivalence between fibrant objects $X\ra Y$
yields a weak equivalence $\Sing X \ra \Sing Y$, and therefore the map
$\Rea^{cof} \Sing X \ra \Rea^{cof} \Sing Y$ is also a weak
equivalence.

For part (ii) recall that any two cosimplicial resolutions of $\gamma$
can be connected by a zig-zag of weak equivalences.  So it suffices to
prove the result in the case where there is a weak equivalence $\Gamma
\ra \Gamma'$.  

We will use $(\re',\sing')$ for the Quillen pair corresponding to
$\Gamma'$.  It is easy to see---using the formulas of \cite[Section
9.5]{D2}, for instance---that there are natural transformations $\re
\ra \re'$ and $\sing' \ra \sing$ induced by $\Gamma \ra \Gamma'$, and
these have the properties that $\Rea A \ra \Rea' A$ and $\Sing' X \ra
\Sing X$ are weak equivalences when $A$ is cofibrant and $X$ is
fibrant.  (In the language of \cite[Def. 5.9]{D2}, this is a {\it
Quillen homotopy\/} from $\re$ to $\re'$.)

If $X$ is fibrant we have a weak equivalence $\Sing' X \ra \Sing X$,
and hence a weak equivalence $Q\Sing' X \ra Q\Sing X$ where $Q$
denotes any cofibrant-replacement functor in $U\cC$. 
Consider the square
\[ \xymatrix{ \Rea Q \Sing' X \ar[r]\ar[d] & \Rea Q \Sing X \ar[d]\\
              \Rea' Q \Sing' X \ar[r] & \Rea' Q \Sing X.
}
\]
All the maps in the square are readily seen to be weak equivalences, and
so we've shown that $\Rea Q \Sing X$ and $\Rea' Q \Sing' X$ are weakly
equivalent via a zig-zag.  By the above proposition, this is what we
wanted.
\end{proof}

\begin{cor}
\label{co:hocosam}
Let $\gamma\colon \cC \ra \cM$ be a functor and $\Gamma\colon \cC \ra
c\cM$ be a cosimplicial resolution of $\gamma$.  Then the induced map
$U\cC \ra \cM$ is homotopically surjective if and only if the natural
maps $\hocolim \ovch \ra X$ are weak equivalences for every fibrant
object $X$.
\end{cor}

\begin{proof}
In light of the above Proposition, this is just a restatement of the
definitions.
\end{proof}

\subsection{Results about canonical homotopy colimits}\fix

In contrast to the canonical homotopy colimit defined above, we can
also consider a more naive construction where we take the
homotopy colimit of the canonical diagram $(\cC \ovcat X) \ra \cM$:
this new object will be denoted \mdfn{$\hocolim (\cC\ovcat X)$}.
Notice the difference between $\hocolim (\cC \ovcat X)$ and $\hocolim
\ovch$---the former was constructed only out of maps $\gamma c \ra X$,
whereas the latter used all `higher-homotopies' $\Gamma^n c \ra X$.
The problem with $\hocolim (\cC \ovcat X)$ is that it is usually not a
homotopy invariant construction---replacing $X$ by another weakly
equivalent object, even if they are both fibrant, may change the
homotopy type of $\hocolim (\cC \ovcat X)$.  On the other hand, this
naive construction is usually much easier to work with than the
canonical homotopy colimit.  There are many instances in this paper
where we get our hands on the canonical homotopy colimit precisely by
showing it agrees with the more naive construction.

Assume now that the image of $\gamma\colon \cC \ra \cM$ is contained
in the cofibrant objects.  In this case we may choose a cosimplicial
resolution $\Gamma$ such that the $0$th object of $\Gamma c$ is
$\gamma c$ itself (rather than an arbitrary cofibrant-replacement).  
There is an obvious map of categories $i\colon(\cC \ovcat X) \ra
(\cC\times \Delta \ovcat X)$ which sends $[c,\gamma c \ra X]$ to the
object $[c\times [0],\Gamma^0 c \ra \gamma c \ra X]$, and this induces
a map of homotopy colimits $i_*\colon\hocolim (\cC\ovcat X) \ra
\hocolim \ovch$.  Our concern will be conditions for this map to be a
weak equivalence.

We have need for one final piece of notation: let $(\cC^n \ovcat X)$
denote the overcategory of $\Gamma^n \colon \cC \ra \cM$: its objects
are pairs $[c,\Gamma^n c \ra X]$.  There is an obvious map $j:(\cC^0
\ovcat X) \ra (\cC^n \ovcat X)$ sending $[c, \Gamma^0 c \ra X]$ to $[c,
\Gamma^n c \ra \Gamma^0 c \ra X]$.  

\begin{prop}
\label{pr:help}
Assume as above that the image of $\gamma\colon \cC \ra \cM$ lies in the
cofibrant objects.  If the maps $j_*\colon \hocolim (\cC^0 \ovcat X)
\ra \hocolim (\cC^n \ovcat X)$ are weak equivalences for all $n$, then
the map $i_*:\hocolim (\cC \ovcat X) \ra \hocolim (\cC\times \Delta
\ovcat X)$ is also a weak equivalence.
\end{prop}

\begin{proof}
Postponed until Section 8.
\end{proof}

The following proposition will be our starting point for obtaining
presentations for combinatorial model categories.  As you will see, it
only concerns the naive construction $\hocolim (\cC\ovcat X)$.  Much
of the work in the rest of the paper involves boot-strapping ourselves
up to a result about canonical homotopy colimits.  (It's useful to
once again compare this result to our discussion in (\ref{se:locdia}).)

\begin{prop}
\label{pr:hoco0}
Let $\cM$ be a combinatorial model category.  Then for sufficiently
large regular cardinals $\lambda$ one has the following: 
\begin{enumerate}[(i)]
\item For every
$X\in \cM$ the canonical map $\hocolim (\cM_\lambda \downarrow X) \ra
X$ is a weak equivalence.
\item The same is true for $\hocolim (\cM_\lambda^{cof} \downarrow X)
\ra X$, where $\cM_\lambda^{cof}$ denotes the full subcategory of
$\cM_\lambda$ consisting of the cofibrant objects.
\end{enumerate}
\end{prop}

\begin{proof}
The proof of $(i)$ is very easy: Since the underlying category of $\cM$ is
locally presentable, for sufficiently large regular
cardinals $\lambda$ the maps $\colim (\cM_\lambda \downarrow X) \ra X$
are isomorphisms for all $X$.  On the other hand the indexing
categories $(\cM_\lambda \downarrow X)$ are $\lambda$-filtered, and in
combinatorial model categories one has that $\lambda$-filtered
colimits are the same as $\lambda$-filtered homotopy colimits for
large enough $\lambda$ (cf. \ref{li:combmc}ii).  So by picking $\lambda$
large enough we may ensure both that the map $\hocolim (\cM_\lambda
\downarrow X) \ra \colim (\cM_\lambda\downarrow X)$ is a weak
equivalence and that the map $\colim (\cM_\lambda\downarrow X) \ra X$
is an isomorphism.  This finishes $(i)$.

For $(ii)$ we must be more careful.  Choose $\lambda$ large enough so
that $(i)$ is satisfied, but also so that $\cM$ has a
cofibrant-replacement functor which maps $\cM_\lambda$ to itself
(\ref{li:combmc}(iii)).  Let $F\colon\cM_\lambda \ra \cM_\lambda$ denote
this functor, let $I=(\cM_\lambda \downarrow X)$, and let
$J=(\cM_\lambda^{cof} \downarrow X)$.

Observe that one has maps $I \llra{f} J \llra{g} I$
where $g$ is the obvious inclusion, and $f$ is the functor sending
$[c,c\ra X]$ to $[Fc,Fc\ra c \ra X]$.  These functors come to us with
natural transformations $gf \ra id$ and $fg \ra id$ induced by the
natural transformation $Fc \ra c$.  Let $D\colon I \ra \cM$
be the canonical diagram sending an object
$[c,c\ra X]$ to $c$.  

The criteria of Proposition~\ref{pr:hocored} are readily checked, and
so we can conclude that $\hocolim_J g^*D \ra\ hocolim_I D$ is a weak
equivalence.  But this is precisely the natural map $\hocolim
(\cM_\lambda^{cof} \downarrow X) \ra \hocolim (\cM_\lambda\downarrow
X)$.  By part ($i$) the codomain is weakly equivalent to $X$, so we
are done.
\end{proof}


\section{The proof for simplicial model categories}

Our goal in this section is to prove Proposition~\ref{pr:hosurj} for
the special case where $\cM$ is a combinatorial model category which
is also a {\it simplicial\/} model category.  The proof for arbitrary
$\cM$ will be given in the next section.  This special case is
presented separately because it is quite a bit easier, yet the steps
are very similar to what we will do for the general case.

\medskip

Choose a regular cardinal $\lambda$ for which
Proposition~\ref{pr:hoco0}(ii) holds: that is, so that the natural map
$\hocolim (\cM_\lambda^{cof} \ovcat X) \ra X$ is a weak equivalence
for all $X$ in $\cM$.  Let $\cC$ denote $\cM^{cof}_\lambda$, for
brevity.  Since $\cM$ is simplicial there is a canonical cosimplicial
resolution for the inclusion $\cC \inc \cM$, and this gives a map
$U\cC \ra \cM$.  The goal will be to show that this map is
homotopically surjective.

\begin{lemma}
\label{le:simp}
The maps $\hocolim (\cC^0\ovcat X) \ra \hocolim (\cC^n\ovcat X)$ are
weak equivalences provided that $X$ is fibrant.
\end{lemma}

\begin{proof}
The objects of $(\cC^n \ovcat X)$ are pairs $[c, c\tens \del{n}\ra X]$
where $c$ is an object of $\cC$ and $c\tens\del{n}\ra X$ is some map
in $\cM$.  Since $\cM$ is simplicial, this map has an adjoint $c \ra
X^{\del{n}}$.  In this way we see that the category $(\cC^n \ovcat X)$
is isomorphic to $(\cC^0 \ovcat X^{\del{n}})$.  The map in which we
are interested is isomorphic to the map
\[ \hocolim (\cC^0\ovcat X) \ra \hocolim (\cC^0\ovcat X^{\del{n}}) \]
induced by $X \ra X^{\del{n}}$.

Now by our choice of $\cC$ we know that $\hocolim (\cC^0\ovcat Z)$ is
naturally weakly equivalent to $Z$, for any $Z$.  So the above map is
weakly equivalent to $X \ra X^{\del{n}}$, which of course is a weak
equivalence because $X$ was fibrant.
\end{proof}

\begin{proof}[Proof of Proposition~\ref{pr:hosurj}, simplicial case]
By Corollary~\ref{co:hocosam}
we must show that for any fibrant $X$ in $\cM$, the natural map
$\hocolim \ovch \ra X$ is a weak equivalence.
Consider the diagram
\[ \xymatrix{ \hocolim (\cC \ovcat X) \ar[r]\ar[dr] 
                    &\hocolim (\cC\times\Delta \ovcat X) \ar[d] \\ 
                    & X.   
}\] 
The above lemma, together with Proposition~\ref{pr:help}, shows that
the horizontal map is a weak equivalence.  The diagonal map is a weak 
equivalence by our choice of $\cC$ (Prop.~\ref{pr:hoco0}). Therefore 
the vertical map is also a weak equivalence, which is what we needed
to prove.
\end{proof}


\section{The proof for non-simplicial model categories}

In this section we prove Proposition~\ref{pr:hosurj} for arbitrary
model categories.  The main surprise is that the category of
`generators' $\cC$ is not chosen to be a subcategory of
$\cM$---instead we have to chose something bigger.

\medskip

\subsection{The plan}
We begin with some general remarks about our approach.  The first
hope would be to take $\cC$ to be the category $\cM_\lambda^{cof}$ for
a sufficiently large cardinal $\lambda$, just as we did for
simplicial model categories.  The difficulty is that we don't know how
to prove Lemma~\ref{le:simp} in this generality.  A second hope might
be to find a method for somehow reducing to the simplicial case, which
we've already handled.

The category of cosimplicial objects $c\cM$ has a natural simplicial
action on it: given an $A \in c\cM$ and a $K\in \sSet$ one can form
new objects $A \tens K$ and $A^K$ (see the appendix).  There is also
an adjoint pair $\ev:c\cM \adjoint \cM:c^*$ where $\ev(A)=A^0$ and
$c^*X$ is the cosimplicial object consisting of $X$ in every
dimension.  In good cases one can find a model structure on $c\cM$ for
which (1) this adjoint pair is a Quillen equivalence, (2) $c\cM$ is a
{\it simplicial\/} model category, and (3) the cofibrant objects of
$c\cM$ are precisely the cosimplicial resolutions.  (This model
structure is dual to the one constructed in \cite{D1}).  Recall that a
cosimplicial resolution is a Reedy cofibrant object of $c\cM$ with the
property that all coface and codegeneracy maps are weak equivalences.

Now if we did have such a model structure on $c\cM$ then we could
apply the result from the previous section to get a presentation for
$c\cM$, and this would also yield a presentation for $\cM$ (using the
Quillen equivalence $c\cM\we \cM$).  Our $\cC$ would be a certain
subcategory of the cofibrant objects in $c\cM$, which are the
cosimplicial resolutions.  Essentially what we do in this section is
unravel this plan in such a way that we never have to actually use the
existence of the model structure on $c\cM$.

\medskip

\subsection{The proof}
Choose a regular cardinal $\lambda$ which is large enough so that
Proposition~\ref{pr:hoco0}(ii) holds, and so that the condition of
Proposition~\ref{pr:combpr}(iii) is satisfied.  Let $\CR$
denote the full subcategory of $c\cM$ consisting of all cosimplicial
resolutions $A^*$ with the property that $A^n \in \cM_\lambda$ for all
$n$.  There is an obvious functor $\gamma\colon \CR \ra \cM$ sending
$A^*$ to $A^0$, and this comes equipped with a natural cosimplicial
resolution $\Gamma\colon \CR \ra c\cM$ which is just the inclusion of
$\CR$ as a subcategory.  These induce a map $U(\CR) \ra \cM$, and the
goal will be to show that this is homotopically surjective.

Let $\Cl$ denote the category $\cM_\lambda^{cof}$ of cofibrant objects
in $\cM_\lambda$, and let $f\colon\CR \ra \Cl$ be the functor sending
$A^*$ to $A^0$.  For any $X\in \cM$ there is an induced map on
overcategories $(\CR\ovcat X) \ra (\Cl\ovcat X)$ sending $[A^*, A^0\ra
X]$ to $[A^0, A^0\ra X]$.

\begin{lemma}
\label{le:nonsimp1}
For any $X\in \cM$, the map $f_*\colon \hocolim (\CR\ovcat X) \ra
\hocolim (\Cl \ovcat X)$ is a weak equivalence.
\end{lemma}

\begin{proof}
Let $\cE$ denote the subcategory of $c\cM$ consisting of the objects
$A^*$ for which $A^n$ is in $\Cl$ for all $n$ (unlike for $\CR$, we
are not requiring $A^*$ to be a cosimplicial resolution).  It is
possible to find a functor $\cR\colon \cE \ra c\cM$ with the following
properties:
\begin{enumerate}[(i)]
\item Each $\cR(A)$ is a Reedy cofibrant object contained in $\cE$;
\item There is a natural weak equivalence $\eta\colon \cR(A) \ra A$;
\item The object of $\cR(A)$ in dimension $0$ is equal to $A^0$, and
the map $\eta\colon\cR(A)^0 \ra A^0$ is the identity.
\end{enumerate}
The map $\cR$ is just a certain Reedy cofibrant-replacement functor
defined on a subcategory of $c\cM$.  In order to construct Reedy
cofibrant-replacements, one starts with the $0$th object and first
makes that cofibrant in $\cM$.  For $A^*\in \cE$ the $0$th object is
already cofibrant, so we can just let it be.  Next one moves
inductively up the cosimplicial object and factors the latching maps
as cofibrations followed by trivial cofibrations (see
\cite[Chap. 5]{Ho}).  By (\ref{pr:combpr}iii) and our choice of
$\lambda$, there are factorization functors which will never take us
outside the category $\cM_\lambda$---this is all that we wanted.

If $A\in \Cl$ then we will write $\cR(A)$ for the result of applying
$\cR$ to the constant cosimplicial object $c^*A$ consisting of $A$ in
every dimension.

Consider the maps $(\Cl\ovcat X) \ra (\CR \ovcat X)$ and $(\CR \ovcat
X) \ra (\Cl \ovcat X)$ induced by $\cR$ and $f$, respectively.  These
maps have the following behavior:
\[ \cR:[c, c\ra X] \mapsto [\cR(c), \cR(c)^0 \llra{id} c \ra X] \qquad 
   f:[A^*, A^0 \ra X] \mapsto [A^0, A^0\ra X].\] 
The composite $f\cR \colon (\Cl\ovcat X) \ra (\Cl \ovcat X)$ is the
identity by property $(iii)$ of $\cR$.  We will show that the other
composite $\cR f$ can be connected to the identity by a zig-zag of
natural transformations, and then we'll apply
Proposition~\ref{pr:hocored}. 

The composite $\cR f$ sends an object $[A^*,A^0 \ra X]$ of $(\CR
\ovcat X)$ to the object $[\cR(A^0),A^0 \ra X]$.  Consider the map
$H\colon (\CR\ovcat X) \ra (\CR\ovcat X)$ which maps $[A^*, A^0 \ra
X]$ to $[\cR A^*,[\cR A]^0 \ra A^0 \ra X]$.  The 
transformation $\eta$ from property $(ii)$ gives a natural
transformation $H \ra \Id$.  On the other hand, for any cosimplicial
object $A^*$ there is a natural map $A^* \ra c^*(A^0)$ and therefore a
map $\cR(A^*) \ra \cR(A^0)$.  This gives a natural transformation
$H\ra \cR f$.  It is easy to check that these transformations
satisfy the conditions of Proposition~\ref{pr:hocored} (see
Remark~\ref{re:hoco} as well).  So we conclude that $f_* \colon
\hocolim (\Cl\ovcat X) \ra \hocolim (\CR\ovcat X)$ is a weak
equivalence, and the same for $\cR_*$ going in the other direction.
\end{proof}

\begin{lemma}
\label{le:nonsimp2}
The canonical map $\hocolim (\CR^0\ovcat X) \ra \hocolim (\CR^n\ovcat
X)$ is a weak equivalence.
\end{lemma}

\begin{proof}
Let $j\colon\del{0}\ra\del{n}$ denote the map of simplicial sets which
includes $\del{0}$ as the last vertex of $\del{n}$.  For any
cosimplicial object $A^*$ there is a corresponding map $j\colon A^0
\ra A^n$; from this we can define a functor $j\colon(\CR^n \ovcat
X) \ra (\CR^0\ovcat X)$ sending the object $[A^*, A^n \ra
X]$ to $[A^*,A^0 \ra A^n \ra X]$.

Let $i\colon (\CR^0 \ovcat X) \ra (\CR^n\ovcat X)$ denote the functor
sending $[A^*, A^0 \ra X]$ to $[A^*, A^n\ra A^0 \ra X]$.  The map we
are concerned with in the statement of the lemma is $i_*$.  Note that
the composite $ji$ is the identity.  We will show that the the other
composite $ij$ can be related to the identity by a zig-zag of natural
transformations, and then we'll apply Proposition~\ref{pr:hocored}.

There is a map of simplicial sets $H\colon\del{n}\times\del{1} \ra
\del{n}$ such that $H$ restricts to the identity map on $\del{n}\times
\{0\}$, and $H$ restricts to the map $j\pi\colon\del{n} \ra\del{0}
\ra\del{n}$ on $\del{n}\times \{1\}$ (left to the reader).  Recall
that if $A^*$ is a cosimplicial object and $K$ is a simplicial set
then one gets a new cosimplicial object $A\tens K$ in a natural way
(see the appendix).  So $H$ induces a map $H\colon
A^* \tens (\del{n}\times \del{1}) \ra A^* \tens \del{n}$, and by
looking at the objects in dimension $0$ we get a map $H\colon (A^*
\tens \del{1})^n \ra A^n$ which is natural in $A^*$.

Consider the functor
$H:(\CR^n\ovcat X) \ra (\CR^n\ovcat X)$ defined by
\[ [A^*, A^n \ra X] \mapsto [A^*\tens \del{1}, (A^*\tens \del{1})^n
\llra{H} A^n \ra X].
\]
The two inclusions $\del{0}\ra\del{1}$ are readily seen to induce
natural transformations $\Id \ra H$ and $ij \ra H$.  The hypotheses of
Proposition~\ref{pr:hocored} are easily checked to hold, and so we may
conclude that $i_*\colon \hocolim (\CR^0\ovcat X) \ra \hocolim
(\CR^n\ovcat X)$, together with $j_*$ going in the other direction,
are both weak equivalences.
\end{proof}

We can now close out the main proof:

\begin{proof}[Proof of Proposition~\ref{pr:hosurj}, general case]
Again, by Corollary~\ref{co:hocosam} we must show that for any fibrant
object $X$ in $\cM$ the natural map $\hocolim (\CR \times \Delta
\ovcat X) \ra X$ is a weak equivalence.

Consider the commutative diagram
\[ \xymatrix{\hocolim (\Cl \ovcat X) \ar[dr]_-{a} 
               & \hocolim (\CR \ovcat X)\ar[l]_-{f_*}\ar[r]^-{i_*}\ar[d] 
               & \hocolim (\CR \times \Delta \ovcat X) \ar[dl]^-{p} \\
 & X.}
\]
Lemma~\ref{le:nonsimp1} says that $f_*$ is a weak equivalence.
Lemma~\ref{le:nonsimp2}, together with Proposition~\ref{pr:help},
implies the same about $i_*$.  Finally, our assumption on $\lambda$
guarantees that $a$ is a weak equivalence (Prop.~\ref{pr:hoco0}).  We
therefore conclude that $p$ is also a weak equivalence, which is what
we wanted.
\end{proof}


\section{More about combinatorial model categories}

\label{ad:combmc}
At this point we have finished with the main ideas of the paper.  In
this section and the next we have only to complete the proofs for some
of the auxiliary results.  This section fills in some of the details
behind the properties of combinatorial model categories singled out in
(\ref{li:combmc}).  The authoritative reference for results like these
will eventually be \cite{Sm}.

\begin{prop}[Smith]
In a combinatorial model category $\cM$ there are functorial
factorizations of a map into a trivial cofibration followed by a
fibration which preserve $\lambda$-filtered colimits for sufficiently
large regular cardinals $\lambda$.  The same is true for the
factorizations as a cofibration followed by a trivial fibration.
\end{prop}

\begin{proof}
The usual factorizations provided by the small object
argument will have the required properties, as long as we use the
transfinite version of the small object argument for a sufficiently
large ordinal.  See \cite{Sm}.
\end{proof}

Proposition~\ref{pr:combpr}(i) is a special case of the above.
Part (iii) of the proposition is the following:

\begin{prop}[Smith]
The factorizations guaranteed by the above proposition have
the following property: for sufficiently large regular cardinals
$\mu$, if $X\ra Y$ is a map between $\mu$-small objects then the
factorizations produce maps $X \trcof \tilde{X} \fib Y$ and $X \cof
\tilde{Y} \trfib Y$ where both $\tilde{X}$ and $\tilde{Y}$ are also
$\mu$-small.   
\end{prop}

In a locally presentable category one can define the {\it size\/} of an
object $X$ to be the smallest regular cardinal $\lambda$ for which $X$
is $\lambda$-small.  The proposition says that past a certain point
the factorizations don't increase size anymore.

\begin{proof}
Pick a regular cardinal $\lambda$ large enough to satisfy the previous
proposition, and also large enough so that $\cM_{\lambda}$ is dense in
$\cM$ (using locally presentability).  The category $\cM_\lambda$ is
small, so applying our given factorizations to maps $X \ra Y$ in
$\cM_\lambda$ only produces a set of new objects.
Therefore there exists a regular cardinal $\nu$ such that applying our
factorizations to maps between $\lambda$-small objects always produces
$\nu$-small objects.

Let $\mu$ be any regular cardinal larger than both $\lambda$ and
$\nu$, and let $X\ra Y$ be a map between $\mu$-small objects.  It
follows from \cite[Prop. 2.3.11]{MP} that we can write $X\ra Y$ as a
colimit of maps $X_\alpha \ra Y_\alpha$ where $X_\alpha$, $Y_\alpha$
are $\lambda$-small and where the indexing category is both
$\mu$-small and $\lambda$-filtered.  Applying our factorization produces
maps $X \trcof \tilde{X} \fib Y$ which are isomorphic to the colimit
of the maps $X_\alpha \trcof \tilde{X}_\alpha \ra Y_\alpha$.  Each
$\tilde{X}_\alpha$ is $\nu$-small (hence $\mu$-small) by our choice of
$\nu$, and so $\tilde{X}$ is a $\mu$-small colimit of $\mu$-small
objects, hence is itself $\mu$-small \cite[Prop. 1.16]{AR}.  This
completes the proof.
\end{proof}

\begin{prop}
Let $\cM$ be a combinatorial model category.  Then for sufficiently
large regular cardinals $\lambda$, $\lambda$-filtered colimits of weak
equivalences are again weak equivalences.
\end{prop}

\begin{proof}
Let $\lambda$ be a regular cardinal large enough so that there are
functorial factorizations preserving $\lambda$-filtered colimits, and
so that the model category has a set of generating cofibrations whose
domains and codomains are $\lambda$-small.  Let $I$ be a
$\lambda$-filtered indexing category, and let $D_1,D_2\colon I \ra
\cM$ be two diagrams.  We suppose that $\eta\colon D_1 \ra D_2$ is a
map of diagrams such that $D_1(i)\ra D_2(i)$ is a weak equivalence for
every $i \in I$, and we'll show that $\colim D_1 \ra \colim D_2$ must
also be a weak equivalence.

Start by factoring the map $\colim D_1 \llra{\eta} \colim D_2$ into a
trivial cofibration followed by a fibration, using our preferred
functorial factorization:
\[ \colim D_1 \trcof X \fib \colim D_2.\]
These maps are the colimits of the maps obtained by applying the
factorization to each spot in the diagram:
\[ D_1(i) \trcof X(i) \fib D_2(i).\]
Now the maps $D_1(i)\ra D_2(i)$ were assumed to be weak equivalences,
and therefore the maps $X(i)\ra D_2(i)$ are actually {\it trivial\/}
fibrations.  If we can show that $\lambda$-filtered colimits of
trivial fibrations are again trivial fibrations then we will be done:
the map $X \ra \colim D_2$ will be a trivial fibration, and so $\colim
D_1 \ra \colim D_2$ will be a weak equivalence.

But we can test if a map is a trivial fibration by checking the
lifting property with respect to our generating cofibrations.  Since
the domains and codomains of these generating cofibrations are
$\lambda$-small, they will factor through some stage of the
$\lambda$-filtered colimit and we will get our lift.
\end{proof}


\section{A leftover proof}
Our goal is to prove Proposition~\ref{pr:help}.  Recall the scene:
$\gamma\colon \cC \ra \cM$ is a functor taking its values in the
cofibrant objects and $\Gamma\colon \cC \ra c\cM$ is a cosimplicial
resolution of $\gamma$ with $\Gamma^0 c=\gamma c$.  There is a
canonical map $\hocolim (\cC\ovcat X) \ra \hocolim \ovch$ for each
$X\in \cM$, and Proposition~\ref{pr:help} gives sufficient conditions
for this map to be a weak equivalence.  The key ingredient in our
proof is the cofibrant-replacement functor $Q$ for $U\cC$, written
down in \cite[Section 2.6]{D2}.

\bigskip

The cosimplicial resolution $\Gamma$ induces a Quillen pair $\re: U\cC
\adjoint \cM: \sing$.  Let $\Sing_n X$ denote the presheaf which forms
the degree $n$ part of $\Sing X$---as usual, we will implicitly
identify $\Sing_n X$ with the corresponding discrete simplicial
presheaf in $U\cC$.

In terms of the above Quillen pair, we know that $\hocolim \ovch$ is
weakly equivalent to $\Rea^{cof} \Sing X$.  We claim that $\hocolim
(\cC\ovcat X)$ is weakly equivalent to $\Rea^{cof} \Sing_0 X$.  The
way to see this is to make use of the functor $Q$ mentioned above.
It's easy to check that $Q \Sing_0 X$ can be identified with the
simplicial replacement of the canonical diagram $(\cC\ovcat X) \ra
\cM$, and then $\Rea Q \Sing_0 X$ gives the usual geometric
realization---the resulting object is precisely $\hocolim (\cC\ovcat X)$.

Moreover, there is an obvious map $\Sing_0 X \ra \Sing X$ obtained by
including the $0$-simplices, and the induced map $\Rea Q \Sing_0 X
\ra \Rea Q\Sing X$ will be weakly equivalent to the map
$\hocolim (\cC\ovcat X) \ra \hocolim \ovch$ we're interested in.

The first thing we will show is that the object $\Rea^{cof} \Sing X$
can be built up as a homotopy colimit of the objects
$\Rea^{cof}\Sing_n X$.  The hypotheses of the Proposition translate
into saying that all the $\Rea^{cof}\Sing_n X$ have the same homotopy
type (as $n$ varies), and so we will be able to collapse the homotopy
colimit down to the $n=0$ piece. 

\bigskip

We start out with a few technical lemmas.  Suppose given a
simplicial diagram $\cF\colon\scop \ra U\cC$.  Since $U\cC$ is a
simplicial model category we may form the geometric realization
$\norm{\cF}$.  On the other hand we may also form the homotopy colimit
using the formula in \cite{BK}, and we will denote this object by
`$\badhocolim \cF$'.  The `bad-' prefix is to remind us that this is
not {\it a priori\/} a homotopy invariant construction, because the
objects in the diagram $\cF$ need not be cofibrant.  Note that there
is a {\it Bousfield-Kan map\/} $\badhocolim \cF \ra \norm{\cF}$, just as
one has in any simplicial category.

\begin{lemma}
For any diagram $\cF$ as above, the Bousfield-Kan map 
\[ \badhocolim \cF \ra \norm{\cF}
\] 
is a weak equivalence in $U\cC$.  
\end{lemma}

\begin{proof}
The point is that the homotopy theory in $U\cC$ all comes from
simplicial sets: the weak equivalences are the objectwise weak
equivalences, and the simplicial structure is the objectwise
structure.  So the lemma is immediately reduced to the corresponding
fact for simplicial sets, which is well-known.
\end{proof}

We will need one other fact about the Bousfield-Kan map in this
context: there are of course natural maps $\cF_0 \ra \norm{\cF}$ and
$\cF_0 \ra \badhocolim \cF$, and it is easy to check that the
triangle
\[ \xymatrix{ \badhocolim\, \cF \ar[r] & \norm{\cF} \\
               \cF_0 \ar[u] \ar[ur]}
\]
is commutative.

\begin{prop}
Let $\cF\colon\scop \ra \cM$ be the diagram given by $[n]\assign \Rea
Q \Sing_n X$.  
There is a commutative triangle
\[ \xymatrix{ \hocolim\, \cF \ar[r]^-\sim & \Rea Q\Sing X \\
               \Rea Q \Sing_0 X \ar[u]^\beta \ar[ur]_\alpha
}
\]
in which the horizontal map is a weak equivalence.
\end{prop}

\begin{proof}
Consider the diagram $\cE\colon\scop \ra U\cC$ given by $[n]\assign Q
\Sing_n X$.  The geometric realization $\norm{\cE}$ is isomorphic to
$Q\Sing X$ (using the definition of $Q$, together with the fact that
for bisimplicial sets the realization is isomorphic to the diagonal).
Our above discussion therefore gives a commutative triangle
\[ \xymatrix{ \badhocolim\, \cE \ar[r] & Q\Sing X \\
               Q \Sing_0 X \ar[u] \ar[ur] 
}
\]
where the horizontal arrow is the Bousfield-Kan map, and therefore a
weak equivalence.  Notice that every object of $\cE$ is cofibrant,
being in the image of $Q$---therefore $\badhocolim \cE$ actually has
the correct homotopy type, and we may drop the `bad-' prefix.
Moreover, $\badhocolim \cE$ is known to be cofibrant in this case.

Now we apply the realization functor to the above triangle, to get
\[ \xymatrix{ \hocolim\, (\Rea Q\Sing_* X) \ar[r]^-\sim & \Rea Q\Sing X \\
               \Rea Q \Sing_0 X \ar[u]^\beta \ar[ur]_\alpha
}
\]
The horizontal map is still a weak equivalence because we applied
$\re$ to a weak equivalence between cofibrant objects.
\end{proof}

\begin{proof}[Proof of Proposition~\ref{pr:help}]
Our assumption is that the maps 
\[ \hocolim (\cC^0\ovcat X) \ra \hocolim (\cC^n\ovcat X) \]
are weak equivalences, for every $n$.  But note that $\hocolim
(\cC^n\ovcat X)$ is precisely $\Rea Q \Sing_n X$ (once again, $Q
\Sing_n X$ may be identified with the simplicial replacement of the
canonical diagram $(\cC^n\ovcat X) \ra \cM$).  The above maps are
weakly equivalent to the iterated degeneracies $\Rea Q \Sing_0 X \ra
\Rea Q \Sing_n X$ in the simplicial object $\Rea Q \Sing_* X$.  From
the fact that these are assumed to be weak equivalences it readily
follows that {\it every\/} map in the simplicial object is a weak
equivalence.  This implies that the natural map $\beta\colon\Rea Q
\Sing_0 X \ra \hocolim [\Rea Q \Sing_* X]$ is also a weak equivalence
(see \cite[Prop. 5.4]{D1}, for instance).

At this point we look at the triangle from the above proposition, and
conclude by the two-out-of-three property that $\alpha\colon\Rea Q
\Sing_0 X \ra \Rea Q \Sing X$ is a weak equivalence.  This is
what we wanted.
\end{proof}


\appendix

\section{Homotopy colimits}

This section has two main goals.  We recall that if $\cM$ is a model
category then the category of cosimplicial objects $c\cM$ has a
natural simplicial structure: if $A\in c\cM$ and $K\in \sSet$ then one
can form objects $A\tens K$ and $A^K$ (in fact all this needs is that
$\cM$ is complete and co-complete).  The tensoring operation is used
in the proof of Lemma~\ref{le:nonsimp2}, and it's also the basis for
the way homotopy colimits are defined in \cite{DHK} and \cite{H}.  We
briefly recall this definition of homotopy colimits, and we list some
basic properties which aren't always stressed.  These properties are
used to prove Proposition~\ref{pr:hocored}, which is a technique for
identifying two homotopy colimits over different indexing categories.
This technique is needed several times in the course of the paper.

\medskip

\subsection{The simplicial structure on \mdfn{$c\cM$}}

If $S$ is a set and $W\in \cM$ then let $W\cdot S$ denote a coproduct
of copies of $W$, one for each element of $S$.  Similarly, let
$W^{\cdot S}$ denote a product of copies of $W$ indexed by the set
$S$.
 
If $K\in \sSet$ and $A\in c\cM$ then we define $A\tens_{\Delta} K \in
\cM$ as the following coend:
\[ A\tens_{\Delta} K = \coeq \Biggl [ \coprod_{[k]\ra[m]} 
A_k\cdot K_m\dbra \coprod_n A_n \cdot K_n \Biggr ]. 
\]
The object $A\tens K \in c\cM$ is then defined to be the
cosimplicial object 
\[ [n] \mapsto A\tens_{\Delta} (K\times \del{n}). \]
The exponential $A^K \in c\cM$ can be defined in a more
straightforward way: it is the cosimplicial object
\[ [n] \mapsto A_n^{\cdot K_n} \]
with the obvious coface and codegeneracy operators.  It is routine
work to check that these definitions give a simplicial structure on
the category $c\cM$---it is exactly dual to the standard simplicial
structure on $s\cM$ (written down in \cite{D1}, for instance).

\subsection{Homotopy colimits}
Suppose that $X\colon I\ra \cM$ is a diagram in a model category, and
let $\Gamma\colon I \ra c\cM$ denote a cosimplicial resolution for
$X$.  The \dfn{homotopy colimit} of $X$ is defined to be the object
\[ \hocolim_I X = \coeq \Biggl [ \coprod_{i \ra j}  \Gamma(i) \tens_{\Delta}
B(j\ovcat I)^{op} \dbra \coprod_{i} \Gamma(i) \tens_{\Delta} B(i\ovcat
\cC)^{op} \Biggr ].
\]
Here $B(i\ovcat I)^{op}$ denotes the classifying space of the category
$(i\ovcat I)^{op}$.  Note that technically speaking the object
$\hocolim_I X$ depends on the cosimplicial resolution $\Gamma$, although
the homotopy type of $\hocolim_I X$ does not.  

This construction of homotopy colimits is essentially the one given in
\cite{DHK} and \cite{H}.  The only difference is that those sources
only require $\Gamma$ to be a cosimplicial {\it framing\/} rather than a
full resolution, which has the effect of giving a construction which
is only homotopy invariant for diagrams of cofibrant objects.  This
distinction is a minor one.

\begin{remark}
The above construction can be seen to have the following properties:
\begin{enumerate}[(i)]
\item If $X_1 \ra X_2$ is a map of $I$-diagrams which is an objectwise
weak equivalence, and $\Gamma_1 \ra \Gamma_2$ is a corresponding map
of resolutions, then the induced map $\hocolim_I X_1 \ra \hocolim_I
X_2$ is also a weak equivalence.
\item Given a functor $f\colon I_1 \ra I_2$ and a diagram $X\colon I_2 \ra
\cM$ one can define a pullback diagram $f^*X\colon I_1 \ra \cM$ by
$f^*X(i)=X(fi)$.  A cosimplicial resolution of $X$ pulls back to a
cosimplicial resolution of $f^*X$, and there is an induced map
  $f_*\colon\hocolim_{I_1} f^*X \ra \hocolim_{I_2} X$.
\item If $I_1 \llra{f} I_2 \llra{g} I_3$ are functors and $X$ is an
$I_3$-diagram (with a chosen cosimplicial resolution) then the
following triangle commutes:
\[ \xymatrix{ \hocolim_{I_1} (gf)^*X \ar[r]^{f_*}\ar[dr]_{(gf)_*} 
                   & \hocolim_{I_2} g^*X \ar[d]^{g_*} \\
                   &\hocolim_{I_3} X.}
\]
\item If $f,g\colon I_1 \ra I_2$ are functors, $X$ is an
$I_2$-diagram, and $\eta:f \ra g$ is a natural transformation, then
$\eta$ induces a map of diagrams $\eta_*\colon f^*X \ra g^*X$.  The
triangle
\[ \xymatrix{ \hocolim_{I_1} f^*X \ar[d]_{\eta_*}\ar[r]^{f_*} 
                   & \hocolim_{I_2} X \\
                   \hocolim_{I_1} g^*X \ar[ur]_{g_*}.}
\]
commutes in the homotopy category.
\end{enumerate}
\end{remark}

The following result and its generalizations (see Remark~\ref{re:hoco}
below) are used several times in the body of the paper.  They are our
main tool for identifying two homotopy colimits over different
indexing categories.

\begin{prop}
\label{pr:hocored}
Let $I$ and $J$ be small categories with a functor $g\colon J \ra I$,
and let $X\colon I \ra \cM$ be a diagram.  Suppose that there is also
a functor $f\colon I \ra J$, together with natural transformations
$\eta\colon gf \ra \Id_I$ and $\theta\colon fg \ra \Id_J$ such that
the following hold:
\begin{enumerate}[(i)]
\item Applying $X$ to the maps $\eta(i)\colon gf(i)\ra i$ yields weak
equivalences, and
\item Applying $X$ to the maps $g(\theta j)\colon g(fg(j)) \ra g(j)$
also yields weak equivalences.
\end{enumerate}
Under these hypotheses, the map $g_*\colon \hocolim_J (g^*X) \ra
\hocolim_I X$ is a weak equivalence.  
\end{prop}

The easiest example to which the result applies is when $f$ and $g$ are
actually an equivalence of categories---in this case conditions (i)
and (ii) are vacuous.  In general the conditions are saying that $f$
and $g$ look like an equivalence as far as $X$ is concerned.

\begin{proof}
Consider the triangle
\[ \xymatrix{ \hocolim_{I} (gf)^*X \ar[d]_{\eta_*}\ar[r]^{f_*}
                  &\hocolim_{J} g^*X \ar[r]^{g_*} 
                   & \hocolim_{I} X \\
              \hocolim_{I} (id)^*X \ar[urr]_{(id)_*} }
\]
which commutes up to homotopy.  The slanted map is the identity, and
the vertical map is a weak equivalence because of assumption (i) on
$X$ (which says that $(gf)^*X \ra X$ is an objectwise weak
equivalence).  So it follows that the composite across the top row 
$(gf)_*$ is a weak equivalence as well.

Likewise, in the triangle
\[ \xymatrix{ \hocolim_{J} (gfg)^*X \ar[d]_{(g\theta)_*}\ar[r]^{g_*}
                  &\hocolim_{I} (gf)^*X \ar[r]^{f_*} 
                   & \hocolim_{J} g^*X \\
              \hocolim_{J} g^*X \ar[urr]_{(id)_*} }
\]
assumption $(ii)$ on $X$ again shows that the composite across the top
is a weak equivalence.  The notation is a little confusing because the
maps labelled $g_*$ in the two diagrams are not exactly the same,
although they are both induced by $g$.  But the maps labelled $f_*$
{\it are\/} the same and this is all we need.  If 
\[ A \llra{a} B \llra{b} C \llra{c} D \]
are maps in some category such that $ba$ and $cb$ are both
isomorphisms, then each of $a$, $b$, and $c$ is also an isomorphism.
Applying this to our situation shows that $f_*$ and the two $g_*$'s
are isomorphisms in the homotopy category of $\cM$, hence they are
weak equivalences.
\end{proof}

\begin{remark}
\label{re:hoco}
In most of the cases where we want to apply this result we actually
don't have a simple natural transformation from the composites $gf$
and $fg$ to the respective identities.  Rather, usually what we have
is a zig-zag of natural transformations.  This is fine, though,
because the same line of reasoning applied to each of the steps in the
zig-zag still shows that the required maps are weak equivalences.
Rather than give a messy formulation of some general result along
these lines, we will leave that to the reader's imagination.
\end{remark}

\bibliographystyle{amsalpha}

\end{document}